\documentclass[a4paper,12pt]{report}
\usepackage[T2A]{fontenc}
\usepackage[cp1251]{inputenc}
\usepackage{amsfonts,amsmath,amssymb}
\usepackage[english]{babel}
\def\beq{\begin{eqnarray}}
\def\eeq{\end{eqnarray}}
\begin{document}
{

\Large{ \bf Bounded solutions of the finite and

\begin{center}    infinite-dimensional dynamical systems.
\end{center}
 }
\small{ 
\begin{center}
Pokutnyi O.A.
\end{center}

Institute of mathematics of NAS of Ukraine, Kiev, Tereshenkivska 3,

E-mail: lenasas@gmail.com

\begin{center}
{\bf Abstract }
\end{center}

Invariant torus are constructed under assumption that the homogeneous system admits an exponential dichotomy on the semi-axes. The main result is closely related with the well-known Palmer's lemma and results of Boichuk A.A., Samoilenko A.M.


{\bf Statement of the problem}\\
Consider the linear inhomogeneous system
\beq
\frac {d\phi}{dt} = a ( \phi ) ,
\qquad
\frac {dx}{dt} = P(\phi) x + f(\phi ),
\label{1}
\eeq
which defined on the direct product of $ m$-dimensional torus $ {\cal T}_m $ or infinite dimensional torus $ {\cal T}_{\infty} $ and the space $R^n$ under assumption  that
$
a ( \phi) \in C^1({\cal T}_m) ;$ $ P(\phi),$ $ f(\phi) \in
C({\cal T}_m);$
 $ \phi = ( \phi_1, \ldots ,\phi_m ) \in
{\cal T}_m ;$  $ x = {\rm col} (x_1, \ldots ,x_n ) \in R^n . $
It is known, that the problem of exsisting and constructing of invariant torus $ x = u(\phi)
\in C({\cal T}_m), $  $ \phi \in {\cal T}_m, $ of the system (\ref{1}) for all $ f(\phi) \in C({\cal T}_m) $ can be solved with using Samoilenko-Green function \cite{S,MSK}. For uniqueness  it is necessary and sufficient for all $ f(\phi) \in C({\cal
T}_m) $ that homogeneous system has no degenerate torus
\beq
\frac {d\phi}{dt} = a (\phi) ,
\qquad
\frac {dx}{dt} = P(\phi) x .
\label{2}
\eeq
It means, that for all $ \phi \in {\cal T}_m $ the system
\beq
\label{3}
\frac{dx}{dt} = P\big (\phi_t(\phi)\big ) x
\eeq
is exponentially-dichotomous
(e-dichotomous) on the whole axis $ R =
(-\infty,+\infty),$ i.e. that there is exists projector $C(\phi) =
C^2(\phi) $ and not dependent from $ \phi,$  $ \tau $ constants $K \geq 1,
\, \alpha > 0 $ such that
\begin{gather}
\begin{gathered}
\displaystyle
\big \Vert\Omega^t_0(\phi)C(\phi)\Omega^0_{\tau}(\phi) \big \Vert
\leq K e^{-\alpha (t-\tau)} ,
\quad t \geq \tau ,
\\
\displaystyle
 \big \Vert \Omega^t_0(\phi)(I-C(\phi))\Omega^0_{\tau}(\phi) \big \Vert
 \leq K e^{-\alpha (\tau -t)},
\quad \tau\geq t ,
\end{gathered}
\label{4}
\end{gather}
for all $ t,\tau \in R ;$  $ \Omega^{t}_{\tau}(\phi)$  $\big (\Omega^{\tau}_{\tau}(\phi) = I_n \big )$ --- is $(n \times n) $-dimensional fundamental matrix of the system~(\ref{3});
$ \phi_t(\phi)-$ is a solution of the Koschi problem
$\dot \phi = a (\phi) ,$  $ \phi_0(\phi) = \phi.$

Consider the case when the system (\ref{3}) doesn't have e-dichotomous on the semi-axes
 $R,$ but e-dichotomous on the semi-axes
$R_+$ and  $R_-$ with projectors $ C_+(\phi)$ and $C_-(\phi) $  $(C^2_{\pm}(\phi) = C_{\pm}(\phi)) $
respectively. It means that \cite{Pal} for the system (\ref{3})
the next inequalities are true
\begin{gather}
\begin{gathered}
\displaystyle
\big \Vert\Omega^t_0(\phi)C_{+}(\phi)\Omega^0_{\tau}(\phi) \big \Vert
\leq K_{1} e^{-\alpha_{1} (t-\tau)},
\quad t \geq \tau  ~~~t, \tau \in \mathbb{R}_{+};
\\
\displaystyle
 \big \Vert \Omega^t_0(\phi)(I-C_{+}(\phi))\Omega^0_{\tau}(\phi) \big \Vert
 \leq K_{1} e^{-\alpha_{1} (\tau -t)},
\quad \tau\geq t, ~~~t, \tau \in \mathbb{R}_{+};
\end{gathered}
\label{400}
\end{gather}


\begin{gather}
\begin{gathered}
\displaystyle
\big \Vert\Omega^t_0(\phi)C_{-}(\phi)\Omega^0_{\tau}(\phi) \big \Vert
\leq K_{2} e^{-\alpha_{2} (t-\tau)} ,
\quad t \geq \tau , ~~~t, \tau \in \mathbb{R}_{-};
\\
\displaystyle
 \big \Vert \Omega^t_0(\phi)(I-C_{-}(\phi))\Omega^0_{\tau}(\phi) \big \Vert
 \leq K_{2} e^{-\alpha_{2} (\tau -t)},
\quad \tau\geq t , ~~~t, \tau \in \mathbb{R}_{-};
\end{gathered}
\label{402}
\end{gather}

is the so called critical case.
In this article it is necessary and sufficient conditions for the existence of invariant torus $ x = u(\phi) \in
C({\cal T}_m),$  $ \phi \in {\cal T}_m ,$ of the system~(\ref{1}) are obtained in that case. It is necessary and sufficient conditions for inhomogeneity $ f(\phi) \in C({\cal T}_m), $ which define invariant manifold are obtained.

{\bf Bounded solutions on the whole axis}

For fixed $\phi
\in {\cal T}_m $ general solutions of the problem
 \beq
\label{5}
\frac{dx}{dt}
 = P\big (\phi_t(\phi)\big ) x + f\big (\phi_t(\phi)\big ),
 \eeq
bounded on the entire semi-axes $R_+$ и $ R_- ,$ have the next form
 \beq
\label{6}
x (t,\phi, \xi_{1})
 = \left \{
\begin {array} {ll}
\displaystyle
\Omega^{t}_{0}(\phi)
C_+(\phi) \xi_{1}
+ \int\limits_ {0}^{t} \Omega^{t}_{\tau}(\phi)
C_+(\phi_{\tau}(\phi)) f (\phi_{\tau}(\phi)) d\tau -
\\[\jot]
\qquad
\displaystyle
-\int\limits_ {t} ^ {\infty} \Omega^{t}_{\tau}(\phi)
(I-C_+(\phi_{\tau}(\phi) ) 
 f (\phi_{\tau}(\phi)) d\tau,
 \quad t \geq 0,
\\[\jot]
\displaystyle
\Omega^{t}_{0}(\phi) (I- C_-(\phi)) \xi_{1}
+ \int\limits_ {-\infty} ^ {t} \Omega^{t}_{\tau}(\phi)
C_-(\phi_{\tau}(\phi)) f (\phi_{\tau}(\phi)) d\tau -
\\[\jot]
\displaystyle
-\int\limits_ {t} ^ {0} \Omega^{t}_{\tau}(\phi)
\big (I-C_-(\phi_{\tau}(\phi) )  f (\phi_{\tau}(\phi)) d\tau,
\quad t \leq 0 ,
\end {array}
\right.
\eeq
for all bounded $f$ and
\beq
\label{600}
x (t,\phi, \xi_{2})
 = \left \{
\begin {array} {ll}
\displaystyle
\Omega^{t}_{0}(\phi)
C_+(\phi) \xi_{2}
+ \int\limits_ {0}^{t} \Omega^{t}_{\tau}(\phi)
(I - C_+(\phi_{\tau}(\phi)) f (\phi_{\tau}(\phi)) d\tau -
\\[\jot]
\qquad
\displaystyle
-\int\limits_ {t} ^ {\infty} \Omega^{t}_{\tau}(\phi)
(C_+(\phi_{\tau}(\phi) ) 
 f (\phi_{\tau}(\phi)) d\tau,
 \quad t \geq 0,
\\[\jot]
\displaystyle
\Omega^{t}_{0}(\phi) (I - C_-(\phi)) \xi_{2}
+ \int\limits_ {-\infty} ^ {t} \Omega^{t}_{\tau}(\phi)
(I - C_-(\phi_{\tau}(\phi)) f (\phi_{\tau}(\phi)) d\tau -
\\[\jot]
\displaystyle
-\int\limits_ {t} ^ {0} \Omega^{t}_{\tau}(\phi)
\big (C_-(\phi_{\tau}(\phi) ) f (\phi_{\tau}(\phi)) d\tau,
\quad t \leq 0 ,
\end {array}
\right.
\eeq
but not for all bounded $f$,
where
 \beq
\label{7}
 C_+(\phi_{\tau}(\phi)) = \Omega^{\tau}_{0}(\phi)C_+(\phi)\Omega^{0}_{\tau}(\phi),
\qquad
C_-(\phi_{\tau}(\phi)) = \Omega^{\tau}_{0}(\phi)
C_-(\phi)\Omega^{0}_{\tau}(\phi) .
\eeq
We say about conditions on $f$ below. Here are some well-known relations \cite{S}
\beq
\begin{array}{c} %
\displaystyle
\Omega _\tau ^t \big ( \phi_s( \phi) \big )
= \Omega _{\tau+s}^{t+s}( \phi),
\qquad
\Omega _{\tau}^{t}( \phi) \Omega ^{\tau }_{s}(\phi)
= \Omega _{s}^{t}( \phi) ,
\\[\jot]
\displaystyle
\big ( \Omega _{\tau }^{t}( \phi) \big ) ^{-1}
= \Omega_{t}^\tau (\phi) ,
\qquad
\phi_{\tau}(\phi_s(\phi)) = \phi_{\tau+s}(\phi) ,
\end{array}
\label{x}
\eeq
which valid for all $ t, \tau, s \in R$, $\phi \in {\cal T}_m$ .

Solutions (\ref{6}) and (\ref{600}) will be bounded on the entire axis $R,$ if the constant vectors
 $ \xi_{1} = \xi_{1}(\phi) \in R^n $ and $ \xi_{2} = \xi_{2}(\phi) \in R^n $
satisfy the next algebraic systems, obtained from the
(\ref{6}) and (\ref{600}) for $t = 0\!:$
\begin{gather}
\displaystyle
\big [C_+(\phi) - (I- C_-(\phi))\big ] \xi_{1}
= \int\limits _ {-\infty}^{0} C_-(\phi)\Omega_{\tau}^{0}(\phi) f
(\phi_{\tau}(\phi)) d\tau +
\nonumber \\
\displaystyle
+ \int\limits _{0}^{\infty}
(I-C_+(\phi))\Omega_{\tau}^{0}(\phi) f (\phi_{\tau}(\phi)) d\tau.
\label{8}
\end{gather}

\begin{gather}
\displaystyle
 \big [C_+(\phi) - (I - C_-(\phi))\big ] \xi_{2}
= \int\limits _ {-\infty}^{0} (I - C_-(\phi))\Omega_{\tau}^{0}(\phi) f
(\phi_{\tau}(\phi)) d\tau +
\nonumber \\
\displaystyle
+ \int\limits _{0}^{\infty}
(C_+(\phi))\Omega_{\tau}^{0}(\phi) f (\phi_{\tau}(\phi)) d\tau.
\label{800}
\end{gather}

Denote by $ D(\phi) = C_+(\phi) -  \big (I- C_-(\phi) \big )$ is $(n\times n)$-dimensional matrix, and by $ D^+(\phi)$~ its Moore-Penrose pseudoinvertible
\cite {Boi}; $ P_{N(D)}(\phi) $
and $ P_{N(D^*)}(\phi) $ are $(n\times n)$-dimensional orthoprojectors:

\begin{gather*}
\displaystyle
P_{N(D)}^2(\phi)
 = P_{N(D)}(\phi) = P_{N(D)}^*(\phi) ,
\\
\displaystyle
P_{N(D^*)}^2(\phi) = P_{N(D^*)}(\phi) = P_{N(D^*)}^*(\phi) ,
\end{gather*}

which project $ R^n $ onto kernel $ N(D) = {\rm ker}\, D(\phi) $ and cokernel $
N(D^*) = {\rm ker}\, D^*(\phi) $ of matrix $ D(\phi) ;$
$$
P_{N(D^*)}(\phi) =
I - D(\phi)D^+(\phi) ,
\qquad  P_{N(D)}(\phi) = I - D^+(\phi)D(\phi).
$$

System (\ref{8}) is solvable if and only if the right part of the system (\ref{8}) belongs to the orthogonal complement of
$ N^{\perp}(D^*(\phi)) = {\rm Im}\, (D(\phi)) $ of the subspace $ N(D^*(\phi)).$
It means that
\begin{gather} %
\displaystyle
P_{N(D^*)}(\phi)   \Bigg \{
\int\limits  _ {-\infty}^{0}
C_-(\phi)\Omega_{\tau}^{0}(\phi) f  \big (\phi_{\tau}(\phi) \big ) d\tau +
\nonumber \\
\displaystyle
+\int\limits  _{0}^{\infty}
(I-C_+(\phi))\Omega_{\tau}^{0}(\phi) f  \big (\phi_{\tau}(\phi) \big ) d\tau
\Bigg\} = 0 .
\label{9}
\end{gather}

In this case the general solutions of the system
(\ref{8}), bounded on the entire axis $R,$ have the  form (\ref{6}) with constant $\xi_{1} = \xi_{1}(\phi) \in R^n, $
which defines from the equation (\ref{8}) by the rule:

\begin{gather} %
\displaystyle
\xi_{1} = D^+(\phi) \Bigg \{
\int\limits  _ {-\infty}^{0}
C_-(\phi)\Omega_{\tau}^{0}(\phi) f (\phi_{\tau}(\phi)) d\tau +
\nonumber \\
\displaystyle
+\int \limits _{0}^{\infty}
 \big (I-C_+(\phi) \big ) \Omega_{\tau}^{0}(\phi) f  \big (\phi_{\tau}(\phi) \big ) d\tau
\Bigg \}
+ P_{N(D)}(\phi) c ,
\quad   c = c(\phi) \in R^n .
 \label{10}
\end{gather}

Substitute (\ref{10}) in (\ref{6}), we obtain that for fixed
$ \phi \in {\cal T}_m $ and inhomogeneity $ f (\phi_{t}(\phi))
\in C({\cal T}_m), $ which satisfies condition (\ref{9}), solutions, bounded on
$ R $ of the system (\ref{5}) have the form

\begin{gather} %
\displaystyle
x (t,\phi, c ) =
\nonumber \\
\displaystyle
=  \Omega^{t}_{0}(\phi)
\begin {cases}
\displaystyle
C_+(\phi) P_{N(D)}(\phi) c + \int\limits  _ {0}^{t}
C_+(\phi) \Omega_{\tau}^{0}(\phi) f (\phi_{\tau}(\phi)) d\tau -
\\[\jot]
\displaystyle
\qquad
-\int\limits  _ {t} ^ {\infty}
 \big (I-C_+(\phi)  \big ) \Omega^{0}_{\tau}(\phi)
 f  \big (\phi_{\tau}(\phi) \big ) d\tau +
\\
\displaystyle
\qquad
+ C_+(\phi)D^+(\phi)
\Bigg \{ \int\limits  _ {-\infty}^{0}
C_-(\phi)\Omega_{\tau}^{0}(\phi) f  \big (\phi_{\tau}(\phi) \big ) d\tau +
\\
\displaystyle
\qquad
+ \int\limits  _{0}^{\infty}  \big (I-C_+(\phi) \big ) \Omega_{\tau}^{0}(\phi) f  \big (\phi_{\tau}(\phi) \big )
d\tau \Bigg \},
\quad t \geq 0,
\\
\displaystyle
 \big (I- C_-(\phi) \big ) P_{N(D)}(\phi) c
 + \int\limits _ {-\infty} ^ {t}
C_-(\phi)\Omega^{0}_{\tau}(\phi) f  \big (\phi_{\tau}(\phi) \big ) d\tau -
\\
\displaystyle
\qquad
-\int\limits  _ {t} ^ {0}
 \big (I-C_-(\phi)  \big ) \Omega^{0}_{\tau}(\phi) f  \big (\phi_{\tau}(\phi) \big )d\tau +
\\
\displaystyle
\qquad
+  \big (I- C_-(\phi) \big ) D^+(\phi)
\Bigg \{ \int\limits  _ {-\infty}^{0}
C_-(\phi)\Omega_{\tau}^{0}(\phi) f  \big (\phi_{\tau}(\phi) \big ) d\tau +
\\
\displaystyle
\qquad
+ \int\limits  _{0}^{\infty}  \big (I-C_+(\phi) \big )
\Omega_{\tau}^{0}(\phi) f  \big (\phi_{\tau}(\phi) \big )
d\tau \Bigg\},
\quad t \leq 0.
\end{cases}
\label{11}
\end{gather}

Since $ P_{N(D^*)}(\phi) D(\phi) = P_{N(D^*)}(\phi)
\big [C_+(\phi)-(I-C_-(\phi))\big ] = 0, $ then
$$
P_{N(D^*)}(\phi)C_+(\phi) = P_{N(D^*)}(\phi)\big (I- C_-(\phi)\big ),
$$
and condition (\ref{9}) is equivalent one of the conditions

\begin{gather}
\begin{gathered}
\displaystyle
P_{N(D^*)}(\phi) \int\limits _ {-\infty}^{+\infty}
C_-(\phi)\Omega_{\tau}^{0}(\phi) f (\phi_{\tau}(\phi)) d\tau = 0,
\\
\displaystyle
P_{N(D^*)}(\phi)\int\limits _{-\infty}^{\infty}
(I-C_+(\phi))\Omega_{\tau}^{0}(\phi) f (\phi_{\tau}(\phi)) d\tau = 0 .
\end{gathered}
 \label{12}
\end{gather}

Since
$$
 \Big [C_+(\phi)-\big (I-C_-(\phi)\big ) \Big ] D^+(\phi)
 = I - P_{N(D^*)}(\phi),
$$
we obtain
$$
C_+(\phi) D^+(\phi) \{ \ldots \} -I \{ \ldots \}
= (I-C_-(\phi)) D^+(\phi) \{ \ldots \},
$$
from the condition (\ref{9}), $\{ \ldots  \}$ is the expression in (\ref{9}).

Since $ D(\phi)P_{N(D)}(\phi) = \big [C_+(\phi)-(I-
C_-(\phi))\big ] P_{N(D)}(\phi) = 0, $ then
$$
C_+(\phi) P_{N(D)}(\phi) = \big (I- C_-(\phi)\big )P_{N(D)}(\phi).
$$
Similarly the system (\ref{800}) is solvable if and only if the right part of system (\ref{800}) belongs to the orthogonal complement of
$ N^{\perp}(D^*(\phi)) = {\rm Im}\, (D(\phi)) $ of the subspace $ N(D^*(\phi)).$
It means that

\begin{gather} %
\displaystyle
P_{N(D^*)}(\phi)   \Bigg \{
\int\limits  _ {-\infty}^{0}
(I - C_-(\phi))\Omega_{\tau}^{0}(\phi) f  \big (\phi_{\tau}(\phi) \big ) d\tau +
\nonumber \\
\displaystyle
+\int\limits  _{0}^{\infty}
(C_+(\phi))\Omega_{\tau}^{0}(\phi) f  \big (\phi_{\tau}(\phi) \big ) d\tau
\Bigg\} = 0 .
\label{900}
\end{gather}

In this case the general solutions of the system
(\ref{800}), bounded on the entire axis $R,$ have the  form (\ref{600}) with constant $\xi_{2} = \xi_{2}(\phi) \in R^n, $
which defines from the equation (\ref{800}) by the rule:

\begin{gather} %
\displaystyle
\xi_{2} = D^+(\phi) \Bigg \{
\int\limits  _ {-\infty}^{0}
(I - C_-(\phi))\Omega_{\tau}^{0}(\phi) f (\phi_{\tau}(\phi)) d\tau +
\nonumber \\
\displaystyle
+\int \limits _{0}^{\infty}
 \big (C_+(\phi))  \Omega_{\tau}^{0}(\phi) f  \big (\phi_{\tau}(\phi) \big ) d\tau
\Bigg \}
+ P_{N(D)}(\phi) c ,
\quad   c = c(\phi) \in R^n .
 \label{1000}
\end{gather}

Substitute (\ref{1000}) in (\ref{600}), we obtain that for fixed
$ \phi \in {\cal T}_m $ and inhomogeneity $ f (\phi_{t}(\phi))
\in C({\cal T}_m), $ which satisfies condition (\ref{900}), solutions, bounded on
$ R $ of the system (\ref{5}) have the form
\pagebreak

\begin{gather} %
\displaystyle
x (t,\phi, c ) =
\nonumber \\
\displaystyle
=  \Omega^{t}_{0}(\phi)
\begin {cases}
\displaystyle
C_+(\phi) P_{N(D)}(\phi) c + \int\limits  _ {0}^{t}
(I - C_+(\phi)) \Omega_{\tau}^{0}(\phi) f (\phi_{\tau}(\phi)) d\tau -
\\[\jot]
\displaystyle
\qquad
-\int\limits  _ {t} ^ {\infty}
 \big (C_+(\phi))  \Omega^{0}_{\tau}(\phi)
 f  \big (\phi_{\tau}(\phi) \big ) d\tau +
\\
\displaystyle
\qquad
+ C_+(\phi)D^+(\phi)
\Bigg \{ \int\limits  _ {-\infty}^{0}
(I - C_-(\phi))\Omega_{\tau}^{0}(\phi) f  \big (\phi_{\tau}(\phi) \big ) d\tau +
\\
\displaystyle
\qquad
+ \int\limits  _{0}^{\infty}  \big (C_+(\phi)  \Omega_{\tau}^{0}(\phi) f  \big (\phi_{\tau}(\phi) \big )
d\tau \Bigg \},
\quad t \geq 0,
\\
\displaystyle
 \big (I - C_-(\phi))  P_{N(D)}(\phi) c
 + \int\limits _ {-\infty} ^ {t}
(I - C_-(\phi))\Omega^{0}_{\tau}(\phi) f  \big (\phi_{\tau}(\phi) \big ) d\tau -
\\
\displaystyle
\qquad
-\int\limits  _ {t} ^ {0}
 \big (C_-(\phi))   \Omega^{0}_{\tau}(\phi) f  \big (\phi_{\tau}(\phi) \big )d\tau +
\\
\displaystyle
\qquad
+  \big (I - C_-(\phi) \big ) D^+(\phi)
\Bigg \{ \int\limits  _ {-\infty}^{0}
(I - C_-(\phi))\Omega_{\tau}^{0}(\phi) f  \big (\phi_{\tau}(\phi) \big ) d\tau +
\\
\displaystyle
\qquad
+ \int\limits  _{0}^{\infty}  \big (C_+(\phi)) \big )
\Omega_{\tau}^{0}(\phi) f  \big (\phi_{\tau}(\phi) \big )
d\tau \Bigg\},
\quad t \leq 0.
\end{cases}
\label{1100}
\end{gather}

Since $ P_{N(D^*)}(\phi) D(\phi) = P_{N(D^*)}(\phi)
\big [I - C_+(\phi) + C_-(\phi))\big ] = 0, $ then
$$
P_{N(D^*)}(\phi)(I - C_+(\phi)) = P_{N(D^*)}(\phi)\big (C_-(\phi)\big ),
$$

and condition (\ref{900}) is equivalent one of the conditions

\begin{gather}
\begin{gathered}
\displaystyle
P_{N(D^*)}(\phi) \int\limits _ {-\infty}^{+\infty}
(I - C_-(\phi))\Omega_{\tau}^{0}(\phi) f (\phi_{\tau}(\phi)) d\tau = 0,
\\
\displaystyle
P_{N(D^*)}(\phi)\int\limits _{-\infty}^{\infty}
(C_+(\phi))\Omega_{\tau}^{0}(\phi) f (\phi_{\tau}(\phi)) d\tau = 0 .
\end{gathered}
 \label{1200}
\end{gather}

Since
$$
 \Big [I - C_+(\phi) + C_-(\phi)\big ) \Big ] D^+(\phi)
 = I - P_{N(D^*)}(\phi),
$$
we obtain
$$
(I - C_+(\phi)) D^+(\phi) \{ \ldots \} -I \{ \ldots \}
= (C_-(\phi)) D^+(\phi) \{ \ldots \},
$$
from the condition (\ref{900}), $\{ \ldots  \}$ is expression in (\ref{900}).

Since $ D(\phi)P_{N(D)}(\phi) = \big [C_+(\phi)-(I-
C_-(\phi))\big ] P_{N(D)}(\phi) = 0, $ then
$$
(I - C_+(\phi)) P_{N(D)}(\phi) = \big (C_-(\phi)\big )P_{N(D)}(\phi).
$$

Consider the case, when homogeneous system (\ref{3}) does not have bounded and unbounded solutions

$$
C_+(\phi) P_{N(D)}(\phi) = \big (I-C_-(\phi)\big ) P_{N(D)}(\phi) = 0.
$$
Then (\ref{11}) and (\ref{1100}) we can rewrite in the form

 \beq
\label{13}
x(t,\phi ) = (G_t(f))(\phi),
 \eeq
\begin{gather*}
\big (G_t(f)\big )(\phi) = \Omega^{t}_{0}(\phi)
\begin{cases}
\displaystyle
\int\limits_ {0}^{t}
C_+(\phi) \Omega_{\tau}^{0}(\phi) f \big (\phi_{\tau}(\phi)\big ) d\tau - \\[\jot]
\displaystyle
\qquad
-\int\limits _ {t} ^ {\infty}
\big (I-C_+(\phi) \big ) \Omega^{0}_{\tau}(\phi)
 f \big (\phi_{\tau}(\phi)\big ) d\tau +
\\[\jot]
\displaystyle
\qquad
+ C_+(\phi)D^+(\phi)  \Bigg \{ \int\limits  _ {-\infty}^{0}
C_-(\phi)\Omega_{\tau}^{0}(\phi) f \big (\phi_{\tau}(\phi)\big ) d\tau +
\\[\jot]
\displaystyle
\qquad
+\int\limits  _{0}^{\infty}\big (I-C_+(\phi)\big )\Omega_{\tau}^{0}(\phi)
f \big (\phi_{\tau}(\phi)\big )
d\tau  \Bigg \},
\quad  t \geq 0,
\\[\jot]
\displaystyle
 \int\limits  _ {-\infty} ^ {t}
C_-(\phi)\Omega^{0}_{\tau}(\phi) f \big (\phi_{\tau}(\phi)\big ) d\tau -
\\[\jot]
\displaystyle
\qquad
- \int\limits  _ {t} ^ {0}
\big (I-C_-(\phi) \big ) \Omega^{0}_{\tau}(\phi)
f \big (\phi_{\tau}(\phi)\big )d\tau +
\\[\jot]
\displaystyle
\qquad
+ \big [ C_+(\phi) D^+(\phi) - I\big ]
 \Bigg \{ \int\limits  _ {-\infty}^{0}
C_-(\phi)\Omega_{\tau}^{0}(\phi) f \big (\phi_{\tau}(\phi)\big ) d\tau +
\\[\jot]
\displaystyle
\qquad
+ \int\limits  _{0}^{\infty}\big (I-C_+(\phi)\big )\Omega_{\tau}^{0}(\phi) f \big (\phi_{\tau}(\phi)\big ) d\tau  \Bigg \},
\quad
t \leq 0 ,
\end{cases}
\end{gather*}

and in the second case

\begin{gather*}
\big (G_t(f)\big )(\phi) = \Omega^{t}_{0}(\phi)
\begin{cases}
\displaystyle
\int\limits_ {0}^{t}
(I - C_+(\phi)) \Omega_{\tau}^{0}(\phi) f \big (\phi_{\tau}(\phi)\big ) d\tau - \\[\jot]
\displaystyle
\qquad
-\int\limits _ {t} ^ {\infty}
\big (C_+(\phi) \big ) \Omega^{0}_{\tau}(\phi)
 f \big (\phi_{\tau}(\phi)\big ) d\tau +
\\[\jot]
\displaystyle
\qquad
+ C_+(\phi)D^+(\phi)  \Bigg \{ \int\limits  _ {-\infty}^{0}
(I - C_-(\phi))\Omega_{\tau}^{0}(\phi) f \big (\phi_{\tau}(\phi)\big ) d\tau 
\\[\jot]
\displaystyle
\qquad
+\int\limits  _{0}^{\infty}\big (C_+(\phi)\big )\Omega_{\tau}^{0}(\phi)
f \big (\phi_{\tau}(\phi)\big )
d\tau  \Bigg \},
\quad  t \geq 0,
\\[\jot]
\displaystyle
 \int\limits  _ {-\infty} ^ {t}
(I - C_-(\phi))\Omega^{0}_{\tau}(\phi) f \big (\phi_{\tau}(\phi)\big ) d\tau -
\\[\jot]
\displaystyle
\qquad
- \int\limits  _ {t} ^ {0}
\big (C_-(\phi) \big ) \Omega^{0}_{\tau}(\phi)
f \big (\phi_{\tau}(\phi)\big )d\tau +
\\[\jot]
\displaystyle
\qquad
+ \big [ (I - C_+(\phi)) D^+(\phi) - I\big ] \\[\jot]
 \Bigg \{ \int\limits  _ {-\infty}^{0}
(I - C_-(\phi))\Omega_{\tau}^{0}(\phi) f \big (\phi_{\tau}(\phi)\big ) d\tau +
\\[\jot]
\displaystyle
\qquad
+ \int\limits  _{0}^{\infty}\big (C_+(\phi)\big )\Omega_{\tau}^{0}(\phi) f \big (\phi_{\tau}(\phi)\big ) d\tau  \Bigg \},
\quad
t \leq 0 ,
\end{cases}
\end{gather*}

will be called generalized Green's operator of the problem about invariant torus of the system (\ref{1}).

Under conditions (\ref{10}), (\ref{1000}) solutions, bounded on $R,$ of the system
(\ref{5}) for fixed $ \phi \in {\cal T}_m $ have the form
(\ref{13}).

We show, that the expression
\beq
\label{19}
x(0,\phi) = u(\phi) = \big (G_0(f)\big ) (\phi) ,
\eeq
which obtained from (\ref{13}) for $ t = 0 , $ define for all $ \phi
\in {\cal T}_m $ invariant torus of the system~(\ref{1}).

{\bf Criterion of existence of invariant torus of nonhomogeneous system}

As shown below, under conditions
 \beq
\label{27}
P_{N(D^*)}(\phi) \int\limits  _ {-\infty}^{+\infty}
C_-(\phi)\Omega_{\tau}^{0}(\phi) f (\phi_{\tau}(\phi)) d\tau = 0,
\eeq

\beq
\label{2700}
P_{N(D^*)}(\phi) \int\limits  _ {-\infty}^{+\infty}
(I - C_-(\phi))\Omega_{\tau}^{0}(\phi) f (\phi_{\tau}(\phi)) d\tau = 0,
\eeq

the nonhomogeneous system (\ref{5}) have bounded solutions on $R$
in the form (\ref{13}) for fixed $\phi \in {\cal T}_m. $
Conditions (\ref{27}) and (\ref{2700}) on solutions $\phi_t(\phi)$
define invariant set. Substitute  $\phi_t(\phi)$ instead of $\phi$ and show that conditions
(\ref{27}) and (\ref{2700}) hold for all $t \in R $ and $ \phi \in
{\cal T}_m. $ From the relations for $ D(\phi) = C_+(\phi)- \big (I-C_-(\phi)\big )$
 we obtain the next equality
\beq \label{20}
 D\big (\phi_t(\phi)\big ) = \Omega^{t}_{0}(\phi) D(\phi)
 \Omega_{t}^{0}(\phi)
 \quad \forall t \in R,
 \quad \forall  \phi \in {\cal T}_m .
 \eeq

Direct check shows, that for all $ t \in R
$ and $ \phi \in {\cal T}_m $ the matrix
 \beq
\label{21}
D^-\big (\phi_t(\phi)\big )
= \big [\Omega^{t}_{0}(\phi) D(\phi)
\Omega_{t}^{0}(\phi)\big ]^-
 = \Omega^{t}_{0}(\phi) D^-(\phi) \Omega_{t}^{0}(\phi)
\eeq

is generalized-invertible to the matrix $ D\big (\phi_t(\phi)\big ) $
and satisfies the next relations \cite{Boi}

\beq
\begin{gathered} %
D^-\big (\phi_t(\phi)\big ) D\big (\phi_t(\phi)\big ) D^-\big (\phi_t(\phi)\big )
= D^-\big (\phi_t(\phi)\big ) ,
\\[\jot]
\displaystyle
D\big (\phi_t(\phi)\big ) D^-\big (\phi_t(\phi)\big )
D\big (\phi_t(\phi)\big ) = D\big (\phi_t(\phi)\big ).
\end{gathered}
\label{211}
 \eeq
From the conditions

\begin{gather*}
\displaystyle
D\big (\phi_t(\phi)\big ) D^-\big (\phi_t(\phi)\big )
= I - P_{N(D)}\big (\phi_t(\phi)\big ),
\\
\displaystyle
D^-\big (\phi_t(\phi)\big ) D\big (\phi_t(\phi)\big )
= I - P_{N(D^*)}\big (\phi_t(\phi)\big )
\end{gather*}

 we get expressions for projectors $ P_{N(D)}(\phi_t(\phi))$ and $P_{N(D^*)}(\phi_t(\phi))$ onto kernel and cokernel of matrix $D(\phi)$ on solutions $ \phi_t(\phi)$ of the respectively Koschi problem for all $ t
\in R$ and $ \phi \in {\cal T}_m \! :$

\beq
\begin{gathered} %
P_{N(D)}\big (\phi_t(\phi)\big )
= \Omega^{t}_{0}(\phi)P_{N(D)}(\phi)\Omega_{t}^{0}(\phi)
=  \Omega^{t}_{0}(\phi)\big [I - D^-(\phi) D(\phi)\big ] \Omega_{t}^{0}(\phi) ,
\\[\jot]
\displaystyle
 P_{N(D^*)}\big (\phi_t(\phi)\big )
 = \Omega^{t}_{0}(\phi)P_{N(D^*)}(\phi)\Omega_{t}^{0}(\phi)
 = \Omega^{t}_{0}(\phi)\big [I - D(\phi) D^-(\phi)\big ] \Omega_{t}^{0}(\phi) .
\end{gathered}
\label{22}
\eeq


We can choose that $ D^-(\phi) = D^+(\phi) $. In that case projectors $P_{N(D)}(\phi) $ and $ P_{N(D^*)}(\phi) $  will be orthoprojectors.

For all  $t \in R $ and $ \phi \in {\cal T}_m $ we have

\begin{gather*}
\displaystyle
P_{N(D^*)}\big (\phi_t(\phi)\big ) \int\limits  _ {-\infty}^{+\infty}
C_-\big (\phi_t(\phi)\big )\Omega_{\tau}^{0}\big (\phi_t(\phi)\big ) f (\phi_{\tau}\big (\phi_t(\phi)\big ))
d\tau =
\\
\displaystyle
 = \Omega^{t}_{0}(\phi)
P_{N(D^*)}(\phi) \int\limits  _ {-\infty}^{+\infty}
C_-(\phi)\Omega_{\tau+t}^{0}(\phi) f (\phi_{\tau+t}(\phi)) d\tau
 = 0,
\end{gather*}

and
\begin{gather*}
\displaystyle
\Omega^{t}_{0}(\phi)
P_{N(D^*)}(\phi) \int\limits  _ {-\infty}^{+\infty}
(I - C_-(\phi))\Omega_{\tau+t}^{0}(\phi) f (\phi_{\tau+t}(\phi)) d\tau
 = 0,
\end{gather*}

From the conditions
(\ref{7}), (\ref{8}), (\ref{22}) we have
$$
u\big (\phi_t(\phi)\big )
= (G_0(f))\big (\phi_t(\phi)\big ) = \big (G_t(f)\big )(\phi)
$$
for all $t \in R $ and $ \phi \in {\cal T}_m. $
It shows that $u\big (\phi_t(\phi)\big ) \in C^1({\cal T}_m) ,$
and the set $u(\phi)$ defines invariant torus of the system (\ref{1}).

In a such way we have the following theorem.


{\bf Theorem.} {\it Let the system \eqref{3} is e-dichotomous on both semi-axes $R_+$ и $ R_-$ with
projectors $C_{\pm}(\phi), $ which satisfy the next equalities for  $\phi_t(\phi) $
$$
C_{\pm}\big (\phi_{t}(\phi)\big )
= \Omega^{t}_{0}(\phi)C_{\pm}(\phi)
\Omega^{0}_{t}(\phi),
\quad  C^2_{\pm}(\phi) = C_{\pm}(\phi).
$$
The system \eqref{1} has invariant torus if and only if nonhomogeneity
$ f(\phi) \in C({\cal T}_m)$ satisfies conditions {\rm \eqref{27}, \eqref{2700}.}
If the homogeneous system {\rm (3)} does not have bounded and unbounded solutions, i.e. the next condition is true
$$
C_+(\phi) P_{N(D)}(\phi) = \big (I-C_-(\phi)\big ) P_{N(D)}(\phi) = 0,
$$
then expression
$$
u(\phi) = \big (G_0(f)\big )(\phi)
$$
which obtain from \eqref{13}, for $ t = 0,$ defines for all $ \phi
\in {\cal T}_m $  invariant torus of system~\eqref{1}.
}

\textbf{\textit{Examples.}}\  \
Consider the problem of the existence of invariant manifold of the system
 \beq
\label{33} \dot
{\varphi} = 1, \\
\\
\qquad \dot {x}_{1}(t) = {\rm th}\, (\varphi)x(t) + f_{1}(\varphi). 
\eeq
\beq
\label{3300} \dot
\qquad \dot {x}_{2}(t) = {\rm - th}\, (\varphi)x_{2}(t) + f_{2}(\varphi).
\eeq
This system has the following characteristics:
$$
\Omega_{0}^{t}(\varphi) = \left(\begin{array}{rcl} \frac{ch(\varphi_{t}(\phi))}{ch\varphi} &
0 \\
0 & \frac{ch\varphi}{ch\varphi_{t}(\varphi)}
\end{array}\right),
$$
$$
C_{+}(\varphi) = \left(\begin{array}{rcl} 0 &
0 \\
0 & 1
\end{array}\right),
C_{-}(\varphi) = \left(\begin{array}{rcl} 1 &
0 \\
0 & 0
\end{array}\right),
$$
$$
D(\varphi) =  \left(\begin{array}{rcl} 0 &
0 \\
0 & 0
\end{array}\right) = D^{+}(\varphi),
$$
$$
P_{N(D)}(\varphi) = P_{N(D^{*})}(\varphi) = I.
$$
In this case the first condition of solvability have the form
$$
\int_{-\infty}^{+\infty} \frac{ch\varphi}{ch(\varphi_{\tau}(\varphi))}f_{1}(\varphi_{\tau}(\varphi))d\tau = 0,
$$
and under this condition the system has the invariant torus in the form
$$
 u_{1}(\varphi) = \left(\begin{array}{rcl} -\int_{0}^{+\infty}\frac{ch(\varphi)f_{1}(\varphi_{\tau}(\varphi))}{ch(\varphi_{\tau}(\varphi))}d\tau \\
0 \end{array}\right),
 $$
and under the second condition of solvability
$$
\int_{-\infty}^{+\infty} ch(\varphi_{\tau}(\varphi))f_{2}(\varphi_{\tau}(\varphi))d\tau = 0
$$
the system has the invariant torus
$$
 u_{2}(\varphi) = \left(\begin{array}{rcl} 0 \\
 -\int_{0}^{+\infty}\frac{ch(\varphi_{\tau}(\varphi))f_{2}(\varphi_{\tau}(\varphi))}{ch(\varphi)}d\tau \end{array}\right).
$$
If, for example $f_{1}(\varphi) = \frac{sh(\varphi)}{ch^{3}(\varphi)}$ and $f_{2}(\varphi) = \frac{sh(\varphi)}{ch^{4}(\varphi)}$, then  the given system has invariant torus in the glued form
$$
u(\varphi) = \left(\begin{array}{rcl}
 -\frac{1}{3ch^{2}(\varphi)} \\
 -\frac{1}{2ch^{3}(\varphi)} \end{array}\right).
$$
We mention that given theory works in the case of infinite dimensional space.
Here is an example.
Consider the countable system of differential equations in the space $BC(\mathbb{R}, l_{2})$ or $BC(\mathbb{R}, l_{2 loc})$ in the next form
$$
\dot{\varphi}(t) = 1,
$$
$$
\frac{dx(t)}{dt} = P(\varphi_{t}(\varphi))x(t) + f(\varphi_{t}(\varphi)),
$$
where
$$
x(t) = (x_{1}(t), x_{2}(t), ...) \in l_{2}~~ \mbox{for all}~~ t,
$$
and
$$
f(\varphi) = (f_{1}(\varphi), f_{2}(\varphi), ...),
$$
where
$$
P(\varphi) = diag \{th(\varphi), th(\varphi), -th(\varphi), -th(\varphi), -th(\varphi), ...\}.
$$
$$
 x(t) = (x_{1}(t), x_{2}(t), ...), f(t) = (f_{1}(t), f_{2}(t), ...)  \in BC(\mathbb{R}, l_{2}).
$$

Here $BC(\mathbb{R}, l_{2}) $ and $BC(\mathbb{R}, l_{2 loc})$ are the spaces of bounded and continuous on the whole axis functions with values in $l_{2}$ or  $l_{2 loc}$.

Matriciant of the system has the form
$$
\Omega_{0}^{t}(\varphi) = diag \{ \frac{ch(\varphi_{t}(\varphi))}{ch\varphi}, \frac{ch(\varphi_{t}(\varphi))}{ch\varphi}, \frac{ch\varphi}{ch(\varphi_{t}(\varphi))}, \frac{ch\varphi}{ch(\varphi_{t}(\varphi))}, \frac{ch\varphi}{ch(\varphi_{t}(\varphi))},... \}.
$$
Projectors have the form
$$
C_{+}(\varphi) = diag\{ 0, 0, 1, 1, ... \}, C_{-}(\varphi) = diag \{1, 1, 0, 0, ...\}.
$$
 Matrixes $D(\varphi) = D^{+}(\varphi) = 0$, and $P_{N(D)} = P_{N(D^{*})} = I$, where $I$ is the identity matrix. Condition of the solvability for the first type of torus has the form
$$
\int_{-\infty}^{+\infty} \frac{f_{i}(\varphi_{\tau}(\varphi))}{ch(\varphi_{\tau}(\varphi))}d\tau = 0,~ i = 1,2,
$$
and invariant torus has the form
$$
x = u_{1}(\varphi) = (- \int_{0}^{+\infty}\frac{ch(\varphi)f_{1}(\varphi_{\tau}(\varphi))}{ch(\varphi_{\tau}(\varphi))}d\tau, -\int_{0}^{+\infty}\frac{ch(\varphi)f_{2}(\varphi_{\tau}(\varphi))}{ch(\varphi_{\tau}(\varphi))}d\tau, 0, ...).
$$
Condition of the solvability for the second type has the form
$$
\int_{-\infty}^{+\infty} f_{i}(\varphi_{\tau}(\varphi))ch(\varphi_{\tau}(\varphi))d\tau = 0,
i \geq 3,
$$
and invariant torus has the form
$$
x = u_{2}(\varphi) = (0, 0, -\int_{0}^{+\infty}\frac{f_{3}(\varphi_{\tau}(\varphi))ch(\varphi_{\tau}(\varphi))}{ch(\varphi)}d\tau, ..., -\int_{0}^{+\infty}\frac{f_{i}(\varphi_{\tau}(\varphi))ch(\varphi_{\tau}(\varphi))}{ch(\varphi)}d\tau, ...),
$$
or, in the glued form 
$$
x = u(\varphi) = (- \int_{0}^{+\infty}\frac{ch(\varphi)f_{1}(\varphi_{\tau}(\varphi))}{ch(\varphi_{\tau}(\varphi))}d\tau, -\int_{0}^{+\infty}\frac{ch(\varphi)f_{2}(\varphi_{\tau}(\varphi))}{ch(\varphi_{\tau}(\varphi))}d\tau, $$
$$ -\int_{0}^{+\infty}\frac{f_{3}(\varphi_{\tau}(\varphi))ch(\varphi_{\tau}(\varphi))}{ch(\varphi)}d\tau, ..., -\int_{0}^{+\infty}\frac{f_{i}(\varphi_{\tau}(\varphi))ch(\varphi_{\tau}(\varphi))}{ch(\varphi)}d\tau, ...).
$$
If, for example,  $f_{i}(\varphi) = \frac{sh\varphi}{ch^{i+2}(\varphi)}$, $i \geq 1$, then we have
$$
x = u_{1}(\varphi) = (-\frac{1}{2ch^{2}(\varphi)}, -\frac{1}{3ch^{3}(\varphi)}, 0, ...),
$$
and
$$
x = u_{2}(\varphi) = (0, 0, -\frac{1}{3ch^{4}(\varphi)},..., -\frac{1}{ich^{i+1}(\varphi)}, ...), $$
or in the glued form 
$$
x = u(\varphi) = (-\frac{1}{2ch^{2}(\varphi)}, -\frac{1}{3ch^{3}(\varphi)}, -\frac{1}{3ch^{4}(\varphi)},..., -\frac{1}{ich^{i+1}(\varphi)}, ...).
$$
Here is denotions as in \cite{Boichuk150}.

\vspace{0.2cm}

{ \small \it Institute of mathematics of NAS of Ukraine, Kiev, 01601, Tereshenkivska str. 3,

 lenasas@gmail.com}


\begin{thebibliography}{99}


\bibitem{S}
{\it Samoilenko A.M.}\ \
Elements of mathematical theory of multifrequency oscillations. -- M.: Science, 1987. -- 304~p.
(in russian).
\bibitem{MSK}
{\it Mitropolsky Yu.O., Samoilenko A.M., Kulik V.L.}\ \
Investigation of dichotomy of linear system of differential equations with Lyapunov
 functions. -- Kiev, 1990. -- 270~p. (in russian)

\bibitem{Pal}
{\it Palmer K. J.}\ \
Exponential dichotomies and transversal homoclinic points //
J. Different. Equat. -- 1984. -- {\bf 55}. -- P.~225\,--\,256.

\bibitem{Bo1}
{\it Boichuk A. A.}\ \
Solutions of weakly nonlinear differential
equations bounded on the whole line // Nonlinear Oscillations. --
 1999. -- {\bf 2}, No~1. -- P.~3\,--\,10.

\bibitem{Boi}
{\it Boichuk A. A., Samoilenko A. M.}\ \
Generalized inverse operators
and fredholm boundary value problems. --
Utrecht; Boston: VSP, 2004. -- 317~p.

\bibitem{B2}
{\it Boichuk A.A. }\ \
Condition of existence of unique Green-Samoilenko function of the invariant torus problem  //
 Ukrainian Math. Journ. --
 2001. -- {\bf 53}, No4. -- p.~556\,--\,559.

\bibitem{B3}
{\it Boichuk A.}\ \
Bounded solutions of differential equations in Banach
space // Colloq. Different. and Difference Equat.
dedicat. Prof. Jaroslav Kurzweil 80-th Birthday: Abstrs (Brno, Czech Republic, Sept. 5\,--\,8, 2006). -- P.~35.

\bibitem{Boichuk150} {\it Boichuk A.A.} Criterion of existence of unique invariant torus of linear extensrion dynamical systems. -- Ukrainian Math. Journal, 2007, V.59, №1. -- p.3 -- 13.


\end{thebibliography}
 \end{document}